\documentclass[11pt]{amsart}
\usepackage{ amsmath, amsthm, amsfonts, hyperref, graphicx, ifpdf}
\usepackage{mathrsfs,epsfig}
\usepackage[dvipsnames,usenames]{color}
\hypersetup{
   unicode=false,          % non-Latin characters in Acrobat bookmarks
   pdftoolbar=true,        % show Acrobat toolbar?
   pdfmenubar=true,        % show Acrobat menu?
   pdffitwindow=false,     % window fit to page when opened
   pdfstartview={},    % fits the width of the page to the window
   pdftitle={},    % title
   pdfauthor={},     % author
   pdfsubject={},   % subject of the document
   pdfcreator={},   % creator of the document
   pdfproducer={}, % producer of the document
   pdfkeywords={}, % list of keywords
   pdfnewwindow=true,      % links in new window
   colorlinks=true,       % false: boxed links; true: colored links
   linkcolor=blue,          % color of internal links
   citecolor=blue,        % color of links to bibliography
   filecolor=blue,      % color of file links
   urlcolor=blue          % color of external links
}

\newtheorem{theorem}{Theorem}[section]
\newtheorem{cor}[theorem]{Corollary}
\newtheorem{lem}[theorem]{Lemma}
\newtheorem{pro}[theorem]{Proposition}
\newtheorem{remark}[theorem]{Remark}
\newtheorem{Def}[theorem]{Definition}
\theoremstyle{definition}

\newcommand{\bnb}{\bar\nabla}

\newcommand{\eus}{Euclidean space}
\newcommand{\ee}{evolution equation}
\newcommand{\es}{evolving surface}
\newcommand{\F}{\frac}
\newcommand{\Fm}{Fuchsian manifold}

\newcommand{\htm}{hyperbolic $3$-manifold}
\newcommand{\hym}{hyperbolic metric}

\newcommand{\kg}{Kleinian group}

\newcommand{\ins}{initial surface}
\newcommand{\mc}{mean curvature}
\newcommand{\mcf}{mean curvature flow}

\newcommand{\nb}{\nabla}

\newcommand{\pc}{principal curvature}

\newcommand{\sff}{second fundamental form}

\newcommand{\tg}{totally geodesic}

\newcommand{\wpp}{warped product}
\newcommand{\wrt}{with respect to}

\newcommand{\hg}{hyperbolic geometry}

\newcommand{\tm}{three-manifold}

\newcommand{\is}{incompressible surface}

\newcommand{\be}{\begin{equation}}
\newcommand{\ene}{\end{equation}}
\newcommand{\br}{\begin{remark}}
\newcommand{\er}{\end{remark}}
\newcommand{\bl}{\begin{lem}}
\newcommand{\el}{\end{lem}}
\newcommand{\bcor}{\begin{cor}}
\newcommand{\ecor}{\end{cor}}
\newcommand{\bpro}{\begin{pro}}
\newcommand{\epro}{\end{pro}}
\newcommand{\ben}{\begin{enumerate}}
\newcommand{\een}{\end{enumerate}}
\newcommand{\bp}{\begin{proof}}
\newcommand{\ep}{\end{proof}}
\newcommand{\bpo}{\begin{pro}}
\newcommand{\epo}{\end{pro}}
\newcommand{\beq}{\begin{equation*}}
\newcommand{\eeq}{\end{equation*}}
\newcommand{\bear}{\begin{eqnarray}}
\newcommand{\eear}{\end{eqnarray}}
\newcommand{\beqar}{\begin{eqnarray*}}
\newcommand{\eeqar}{\end{eqnarray*}}
\newcommand{\bt}{\begin{theorem}}
\newcommand{\et}{\end{theorem}}

\newcommand{\nablabar}{{\overline{\nabla}}}

\newcommand{\n}{\mathbf{n}}
\newcommand{\vnu}{\boldsymbol{\nu}}

\newcommand{\R}{\mathbb{R}}
\newcommand{\I}{\mathbb{I}}
\renewcommand{\H}{\mathbb{H}}

\newcommand{\inner}[2]{\langle #1\,,#2\rangle}
\newcommand{\ddl}[2]{\frac{d{#1}}{d{#2}}}

\newcommand{\ppl}[2]{\frac{\partial{#1}}{\partial{#2}}}

\numberwithin{equation}{section}

\allowdisplaybreaks

\def\XXint#1#2#3{{\setbox0=\hbox{$#1{#2#3}{\int}$}
    \vcenter{\hbox{$#2#3$}}\kern-.5\wd0}}

\makeatletter
\def\@citestyle{\m@th\upshape\mdseries}
\def\citeform#1{{\bfseries#1}}
\def\@cite#1#2{{%
  \@citestyle[\citeform{#1}\if@tempswa, #2\fi]}}
\@ifundefined{cite }{%
  \expandafter\let\csname cite \endcsname\cite
  \edef\cite{\@nx\protect\@xp\@nx\csname cite \endcsname}%
}{}
\makeatother
\def\<{\langle}
\def\>{\rangle}
\def\({\left(}
\def\){\right)}

\def\Ric{{\rm Ric}}

\def\p{\partial}

\begin{document}

\title[MCF in Fuchsian manifolds]{Mean curvature flow in Fuchsian manifolds}

\author{Zheng Huang}
\address[Z. ~H.]{Department of Mathematics, The City University of New York, Staten Island, NY 10314, USA}
\address{The Graduate Center, The City University of New York, 365 Fifth Ave., New York, NY 10016, USA}
\email{zheng.huang@csi.cuny.edu}

\author{Longzhi Lin}
\address[L.~L.]{Mathematics Department\\University of California, Santa Cruz\\1156 High Street\\
Santa Cruz, CA 95064\\USA}
\email{lzlin@ucsc.edu}

\author{Zhou Zhang}
\address[Z. ~Z.]{School of Mathematics and Statistics, The University of Sydney, NSW 2006, Australia}
\email{zhangou@maths.usyd.edu.au}

\date{\today}
\subjclass[2010]{Primary 53C44, 57M10}

%------------------------------------------------------

\begin{abstract}

Motivated by the goal of detecting minimal surfaces in hyperbolic manifolds, we study geometric flows in complete hyperbolic $3$-manifolds. In general, the flows might develop singularities at some finite time. In this paper, we investigate the {\mcf} in a class of complete {\htm}s ({\Fm}s) which are warped products of a closed surface of genus at least two and $\R$. We show that for a large class of closed initial surfaces, which are graphs over the {\tg} surface $\Sigma$, the {\mcf} exists for all time and converges to $\Sigma$. This is among the first examples of converging {\mcf}s starting from closed hypersurfaces in Riemannian manifolds. We also provide calculations for the general warped product setting which will be useful for further works.

\end{abstract}

\maketitle

%------------------------------------------------------

\section {Introduction}
\subsection{The setting}

The {\mcf} has been studied extensively in various ambient Riemannian manifolds and in most cases the flow of closed submanifolds develops singularities in finite time by the avoidance principle for the mean curvature flow. For instance, we know that the finite time singularity has to occur for any compact initial hypersurface in {\eus} under the {\mcf} (\cite{Hui84}). The study of singularity formation has been a focal point of the field, see for instance \cite{Hui86, HS09, CM12} and many others. Hyperbolic manifolds are known to possess extremely rich geometric structures. Our motivation is to better understand the {\mcf} in hyperbolic manifolds, hoping that in further studies we can detect interesting geometric objects by running the mean curvature flow or similar flows to time infinity without developing any singularity or after handling possible singularities.

As a first step, we focus on the {\Fm}s in this paper. {\Fm}s are probably the most elementary complete, non-simply 
connected {\htm}s. A {\Fm} is obtained as a quotient space of $\H^3$ by a Fuchsian group. Let $M^3$ be a {\Fm}, 
and we always assume the genus of any {\is} of $M^3$ is at least two so that it carries its own {\hym}. From differential
geometry point of view,  it is a warped product of a hyperbolic surface $\Sigma$ with $\R$, with the metric
$$ds^2 = dr^2 + \cosh^2(r)g_0\,,$$ where $g_0$ is the {\hym} on $\Sigma$. Therefore the surface $\Sigma$
is {\tg} in $M^3$. Clearly it is the only such surface in $M^3$.

Our main analytical tool is the {\mcf} equation, which has the following form:
\be\label{mcf}	
   \left\{
   \begin{aligned}
      \ppl{}{t}\,F(x,t)&=-H(x,t)\vnu(x,t)\ ,\\
      F(\cdot,0)&=F_{0}\ ,
   \end{aligned}
   \right.
\ene
where $H(x,t)$ and $\vnu(x,t)$ are the {\mc} and unit normal vector respectively at $F(x,t)$ of the evolving surface $S(t)$, and our convention of the {\mc} is the sum of the {\pc}s. 

\begin{Def}
A smooth closed surface $S_0$ in $M^3$ is a graph over the {\tg} surface $\Sigma$ if there is a constant $c_0 > 0$ such that the angle function $\Theta_0 = \langle{\n},{\vnu_0}\rangle \geq c_0$, where $\langle\cdot, \cdot\rangle$ is the metric over $M^3$, $\n = \ppl{}{r}$ is the unit vector field over $M$ which is perpendicular to $\Sigma$ and $\vnu_0$ is the unit normal vector on $S_{0}$ of our choice. 
\end{Def}
Note that $\Theta_0 \in (0,1]$ if $S_0$ is a graph, and $\Theta_0 \equiv 1$ if and only if $S_0$ is \emph{equidistant} from $\Sigma$ (sometimes called \emph{parallel} to $\Sigma$, or a \emph{level surface} to $\Sigma$).

The {\mcf} in {\wpp} manifolds was also investigated by other authors, see e.g. \cite{BM12}. Note that the 
{\wpp} structures in \cite{BM12} are completely different from ours. Geometrically, their warped structure can be thought as a real line bundle over a surface, while ours is a surface bundle over the real line. So the evolving hypersurfaces in their case are \textit{equidistant graphs} over a reference hypersurface, while our evolving hypersurfaces are the more natural \textit{geodesic graphs} (i.e. graphs over totally geodesic surfaces). One does not in general expect a geodesic graph to stay geodesic graphs under the {\mcf}. In fact, it was mentioned in \cite{BM12} that in \cite{Unt98} Unterberger gave an example of hypersurface which is a geodesic graph but loses this graphical property when it evolves under the mean curvature flow. On the other hand, we will see in this work that there indeed exists a large class of closed initial surfaces $S_0$'s in Fuchsian manifolds, as geodesic graphs over the {\tg} surface $\Sigma$ with an explicit lower bound on the angle of $S_0$, such that the {\mcf} starting from $S_0$ remains as geodesic graphs for all time and converges smoothly to $\Sigma$.

%%%%%%%%%%%%%%%%%%%%%%%
\subsection{Main Result}
In this paper, we prove that if the angle function on the initial surface has a positive lower bound depending only 
on its maximal distance to the reference surface $\Sigma$, then the {\mcf} with such an {\ins} exists for all time 
and converges to the {\tg} surface $\Sigma$ in a {\Fm} 
$$(M^3, g_M) =  \bigl(\R\times_{\cosh(r)}\Sigma, dr^{2}+\cosh^{2}(r)g_0\bigr).$$ 
More precisely, we have
\bt\label{main}
Let $M^3$ be a {\Fm} and $\Sigma$ the unique closed {\tg} surface in $M^3$. Then for any $a_0>0$, if the initial smooth closed surface $S_0 \subset M^3$ has hyperbolic distance no larger than $a_0$ to $\Sigma$ and the minimum of the initial angle satisfies 
\be \label{assumption}
\min_{p\in S_0} \Theta_0(p) \ge \tanh(a_0),
\ene
then the {\mcf} with initial surface $S_0$ exists for all time, remains as geodesic graph over $\Sigma$ and converges continuously to $\Sigma$. Moreover, the convergence is smooth if the above inequality is strict.
\et

\br Notably, in the special case where $\min_{p\in S_0}\Theta_0=1$, namely the {\ins} $S_0$ is equidistant from the {\tg} 
surface $\Sigma$, the evolving surface $S(t)$ remains equidistant from $\Sigma$, i.e. $\Theta(t) = \langle{\n},{\vnu}\rangle (t) \equiv 1$ for all $t \ge 0$ and converges to the totally geodesic surface $\Sigma$ by the explicit solution \eqref{ode}, which is the prototype motivating our consideration. The result itself illustrates the interaction between geometric data of the flow and the ambient space. Namely, the lower bound of the angle function $\Theta_0$ of $S_0$ on the right hand side of \eqref{assumption} is nothing but the principle curvature of the equidistant surface $\Sigma(a_0)$, where $a_0$ is the maximal distance of $S_0$ to $\Sigma$. At this moment, it is more like a coincidence as the optimal choice by our argument. We provide some discussion regarding the formation of singularities for the graph case in general in Section \ref{sec4} and hope to sort out the underlying connection in future works. 
\er

Our techniques can be generalized to higher dimensional {\wpp} manifolds of similar structure with appropriate variations of curvature conditions. 
%%%%%%%%%%%%%%%%%%%%%%%%%%%%%
\subsection{Interaction between geometry and analysis}

We want to highlight the interaction between analytical methods and geometric structures. Our setting of {\Fm}s 
allows us to take advantage of its {\hg} in several stages of this work. It is a basic fact that the level sets 
$\{(\Sigma(r), \cosh^2(r)g_0)\}_{r\in\R}$ of the {\tg} surface $\Sigma = \Sigma(0)$ form an equidistant foliation 
of the {\Fm} $M^3$. Moreover, each fiber of the foliation is umbilic, which enables us to obtain the {\mcf} with {\ins} 
$\Sigma(r)$ an explicit solution for any fixed $r \in \R$ (see \eqref{ode}). We use this special {\mcf} as 
barriers and the avoidance principle for the {\mcf}s (see e.g. \cite{Hui86}) to push flow to the destination 
$\Sigma$ (see \S 3.1). Furthermore, we use the presence of a special vector field $V = \cosh(r)\ppl{}{r}$ 
(see \eqref{killing}) in a {\Fm} to derive explicitly the {\ee} for the angle $\Theta$ (see \eqref{ee-theta}).
%%%%%%%%%%%%%%%%%%%%%%%

%------------------------------------------------------
\subsection{Outline of the paper}
We provide some preliminary results in \S2. Heart of the matter is to prove the preservation of graphical 
property of the flow, and we prove our main result Theorem ~\ref{main} in \S 3. The scheme is the following. We first show the evolving surface must stay in a bounded region in $M^3$ for all time (the Squeeze Lemma ~\ref{bd-u}) as long as the flow exists, then we derive the {\ee} for the angle function $\Theta(\cdot,t)$ (Theorem ~\ref{ee-theta}), and then most of the work is devoted to prove that the {\es}s stay graphical under the initial condition on distance and angle. Note that the uniform positive lower bound of the angle function $\Theta = \langle{\n},{\vnu}\rangle$ locally gives the uniform $C^1$-estimate of the graph function which represents the evolving surface. Therefore once we have established uniform bound for $\Theta$, standard parabolic theory (\cite {LSU68}) gives bounds for all higher derivatives. In particular the {\sff} for the {\es} $S_t$ in $M^3$ is uniformly bounded. Huisken's theorem (\cite{Hui86}) then guarantees that the {\mcf} exists for all time. The Squeeze Lemma ~\ref{bd-u} then gives the convergence of the flow. The smoothness will follow.

Finally, in \S 4 we remove the assumption of the angle function on the initial graph and illustrate the possible formation of singularities.

%------------------------------------------------------
\subsection{Acknowledgements} We are grateful to Hengyu Zhou and Biao Wang for many useful discussions and helpful suggestions. We also thank the referee for insightful comments which help to clarify several statements in the paper. Z. Huang acknowledges support from U.S. NSF grants DMS 1107452, 1107263, 1107367 ``RNMS: Geometric Structures and Representation varieties" (the GEAR Network) and a grant from the Simons Foundation (\#359635). L. Lin was partially supported by a Faculty Research Grant awarded by the Committee on Research from UC, Santa Cruz. Z. Zhang was partially supported by Australian Research Council Discovery Project DP110102654 and Future Fellowship FT150100341. This work was partially supported by the National Science Foundation under Grant No. DMS-1440140 while the authors were in residence at the Mathematical Sciences Research Institute in Berkeley, California, during the Spring 2016 semester.
%------------------------------------------------------

\section{Preliminaries}\label{Preliminaries}
In this section, we fix our notations, and introduce some preliminary facts that will be used throughout.

%------------------------------------------------------

\subsection{Fuchsian manifolds}

A {\Fm} $M^3$ is defined as a {\wpp} space $\R \times \Sigma$. This is a complete {\htm} of fundamental
importance in {\hg} and {\kg} theory. The metric on a {\Fm} $M^3$ is explicitly given as:
\be\label{metric}
(M^{3}, g_M) = (\R\times_{\cosh(r)}\Sigma, dr^{2}+\cosh^{2}(r)g_0),
\ene
where $g_0$ is the induced metric on the surface $\Sigma$ which is hyperbolic. A fascinating feature of
the geometry of the {\Fm} $M^3$ is that the equidistant surface $\Sigma(r)$ of the {\tg} surface $\Sigma(0) = \Sigma$
forms a global foliation of $M^3$, and each fiber surface $\Sigma(r)$ is umbilic, with constant {\pc} of 
$\tanh (r)$, cf. \cite{One83}.

Another important fact which will be useful for us is the existence of the special vector field $V$ in $M^3$ (see for instance \cite{Bre13}). Namely, 
\be \label{killing1}
V = \cosh(r)\n
\ene
satisfies that
\be\label{killing}
\bnb_{X}V = \sinh(r)X,
\ene
for any smooth vector field $X$, and $\bar\nabla$ is the Levi-Civita covariant derivative in $M^3$. In the paper, notations with a overhead ``-" stand for quantities for the ambient manifold $M^3$. Making use of this $V$, one can deduce the following well-known formula. 
\bl\label{generalid}
Let $X$ be any tangent vector field in a {\Fm} $M^3$, then we have
\be\label{id}
\bnb_X\n = \tanh(r)(X-\langle X,\n\rangle\n).
\ene
\el
\bp
Since both sides of the identity are linear in $X$, we only have to verify the identity by taking a local frame.
Fix any point $\bar{q} \in M^3$ and let $\{\bar{e}_1, \bar{e}_2\}$ be a local normal frame at the corresponding 
point $q \in \Sigma$, then $\{ \bar{e}_1, \bar{e}_2, \bar{e}_3 = \n\}$ forms a local normal frame at $\bar{q}$. Then 
\eqref{id} is obviously true for $X = \n$. When $X = \bar{e}_1$, we have
\beq
\bnb_{\bar{e}_1}\n = \bar\Gamma_{31}^3\n+ \bar\Gamma_{31}^j\bar{e}_j.
\eeq
Since the metric on $M^3$ is given explicitly as \eqref{metric}, its Christoffel symbols $\bar\Gamma_{ij}^k$
can be computed explicitly. In our case, we have $\bar\Gamma_{31}^3 = 0$ and
$\bar\Gamma_{31}^j\bar{e}_j = \tanh(r)\bar{e}_j$. One sees that \eqref{id} holds for $X = \bar{e}_1$. 
Similarly we can verify the case of $X = \bar{e}_2$.
\ep

%%%%%%%%%%%%%%%%%%%%%%%%%%%%%%%%%%%%%%%%%%%%
\subsection{Mean curvature flow}\label{mcf-s}
Let $F_{0}: S \to{}M^3$ be the immersion of an {\is} $S$ in a {\Fm} $M^3$. We
assume that $S_{0}=F_{0}(S)$ is a graph over $\Sigma$ with respect to $\n$, i.e. $\langle{\n},{\vnu_0}\rangle \geq c_0>0$, here
$c_0$ is a constant to be determined later.

We consider a family of immersions of surfaces in $M^3$ moving under the mean curvature flow \eqref{mcf}:
\beq
   F:S\times[0,T)\to{}M^3\ , \quad{}0\leq{}T\leq\infty\\
\eeq
with $\ppl{}{t}\,F(x,t)=-H(x,t)\vnu(x,t)$, and $F(\cdot,0)=F_{0}$. For each $t\in[0,T)$, $S(t)=\{F(x,t)\in{}M^3\ |\ x\in{}S\}$ is the evolving 
surface at time $t$, and $S(0)$ is the initial surface $S_0$.

The short-time existence of the solutions to \eqref{mcf} is standard for closed immersions, see e.g. \cite{HP96}. 
Initial compact surface develops singularities in finite time along the mean curvature flow in {\eus}, and in fact 
the norm of the {\sff} blows up if the singularity occurs in finite time, see \cite {Hui84, Hui86}.
%------------------------------------------------------
%------------------------------------------------------

\subsection{Evolution equations}
In this subsection, for completeness we collect and derive a number of {\ee}s of some quantities on
$S(t)$, $t\in[0,T)$, which are involved in our calculations. We include here standard evolution equations for the 
{\mc} $H(\cdot,t)$, and the square norm of the {\sff} $|A(\cdot,t)|^2$.
\bt Along the \mcf, one has 
\begin{align}
      \left(\ppl{}{t}-\Delta\right)\,H=
        &\,H(|A|^{2}-2)\ ,\\
            \left(\ppl{}{t}-\Delta\right)\,|A|^{2}=&\, -2|\nabla{}A|^{2}+2|A|^4+4(|A|^2 -H^2)\ .
\end{align}
\et
\bp
These equations are deduced for general Riemannian manifolds in \cite{Hui86}. In our case of hyperbolic {\tm},
the ambient space $M^3$ has constant sectional curvature $-1$, and the Ricci curvature $Ric(\vnu,\vnu) = -2$ for
any unit normal vector $\vnu$.

The lemma then follows from combining these explicit curvature conditions and curvature equations
$\bar{R}_{3i3j} = -g_{ij}$, as well as the well-known Simons Identity (see e.g. \cite{Sim68} or \cite{SSY75}), satisfied
by the {\sff} $a_{ij}$.
\ep

Our estimates are mainly for the {\it height function} $u(x,t)$ and the {\it angle function}, which is the cosine of the geometric angle, $\Theta(x,t)$ on $S(t)$:
\begin{align}
   u(x,t)&=r (F(x,t)), \\
      \label{eq:gradient function}
   \Theta(x,t)&=\langle \vnu, \n \rangle (F(x,t)),
\end{align}
where $r(p) = \pm {\text{dist}}(p,\Sigma)$ for all $p \in M^3$, the {\it signed distance} to the fixed reference surface $\Sigma$.

We have $\Theta (x,t) \in [0,1]$ by the choice. It is clear that $S(t)$ is a geodesic graph over $\Sigma$ if $\Theta > 0$ on $S(t)$. $\Theta(\cdot, 0)=\Theta_0$ is for the initial surface.

\bt
\label{old-eq}
The {\ee}s of $u$ and $\Theta$ have the following form:
\begin{align} \label{eq-u0}
   \ppl{}{t}\,u
        =&\,-H\Theta,\\
               \label{ev-theta}
    \left(\ppl{}{t}-\Delta\right)\Theta=
       &\,(|A|^{2}-2)\Theta+\n(H_{\n})
          -H\inner{\nablabar_{\vnu}\n}{\vnu}\\ 
           \label{ev-theta-1}        
        =&\, (|A|^2-2)\Theta+\frac{1}{2}(\bar\nabla_{\vnu} L_{\n} g)(e_i, e_i) -(\bar\nabla_{e_i}L_{\n} g)(\vnu, e_i)\\
        &\, -a_{ij}L_{\n} g(e_i, e_j), 
          \notag% \label{ev-theta-2}        
\end{align}
where $\n(H_{\n})$ is the variation of {\mc} of $S(t)$ under the deformation vector field $\n$ and $L_{\n} g$ is the Lie derivative of the metric $g$ in the direction $\n$.
\et
\bp
The first equation is self-evident. The second equation can be found in (\cite {Bar84}, \cite{EH91}) in the Lorentzian setting. We include the proof for the Riemannian setting in the Appendix for the sake of completeness. 
\ep
\br
As we will see in the Appendix, the equation for $\Theta$ in Theorem \ref{old-eq} is quite general, and thus very difficult to work with, 
especially the term $\n(H_{\n})$. One of the key observations in our work is that we can take advantage of the explicit nature of the 
ambient {\wpp} metric and derive more workable equations in our case (see \eqref{evo-Theta}).
\er
%------------------------------------------------------

\section{Proof of Main Theorem}\label{sec3}
We prove the main theorem in this section. In \S3.1, we use the explicit solution to the {\mcf} when the {\ins} is 
parallel to the {\tg} surface $\Sigma$ to conclude (Theorem ~\ref{bd-u}) that the {\mcf} starting from any 
closed {\ins} $S_0$ stays in a compact region in $M^3$ as long as it exists. This is standard $C^0$-estimate using the avoidance principle for the {\mcf}.  In \S3.2, we derive the key {\ee} 
(Theorem ~\ref{ee-theta}) for the angle function $\Theta(\cdot,t)$. In \S 3.3, we prove the preservation of graphical 
property. We finally prove Theorem ~\ref{main} in \S 3.4.
%%%%%%%%%%%%%%%%%%%%%%%%%%
\subsection{The squeeze}
The following theorem is probably known previously, but we include the proof here for the sake of completeness.
\bt\label{bd-u}
Let $S_0$ be as in Theorem \ref{main}, then as long as the flow exists we have
\be\label{expdecay}
-\sinh^{-1}(e^{-2t}\sinh(a_0)) \le u(\cdot,t) \leq \sinh^{-1}(e^{-2t}\sinh(a_0)).
\ene
\et
\bp
It follows from basic {\hg} that, if we denote $\Sigma(r)$ (resp. $\Sigma(-r)$) the parallel surface equidistant 
$r$ to $\Sigma(0) = \Sigma$ which stays in the positive (resp. negative) side of $\Sigma$, then $\Sigma(r)$ 
(resp. $\Sigma(-r)$) is an umbilic surface of constant {\pc} $\tanh(r)$ (resp. $-\tanh(r)$).

Now we consider the  {\mcf} with {\ins} $\Sigma(a_0)$ such that $a_0 \geq \text{dist} (x,\Sigma)$ for any 
$x \in S_0$. By the well known uniqueness of the flow, the {\mcf} equation is reduced to the following ODE of $R(t)$, the $r$-value of the evolving equidistant surface:
\beq
\frac{dR}{dt} = -2\tanh(R),
\eeq
with initial condition $R(0) = a_0$, which yields an explicit solution:
\be\label{ode}
R(t) = \sinh^{-1}(e^{-2t}\sinh(a_0)).
\ene
Similar calculations hold for $\Sigma(-a_0)$. One sees that such {\mcf} exists for all time, and all {\es}s are 
umbilic and converge to $\Sigma$ as $t \to \infty$.

Now by assumption, the {\ins} $S_0$ lies between umbilic slices $\Sigma(a_0)$ and $\Sigma(-a_0)$, the 
conclusion then follows from the avoidance principle for the {\mcf}.
\ep

Next we derive a general equation for $\Delta u$.
\bl\label{equation-u}
Let $S \subset M^3$ be a closed surface that is a geodesic graph over $\Sigma$, and $u(x)$ is the signed 
distance of $x \in S$ to $\Sigma$. Then we have:
\be\label{eq-u1}
\Delta u = \tanh(u)(1+\Theta^2) - H\Theta\,,
\ene
where $\Delta$ is the Laplace operator on $S$ {\wrt} the induced metric.
\el
\bp
For any point $x \in S$, choose $\{e_1, e_2\}$ (with $e_3= \vnu$) to be a local normal frame of $S$ at $x$.  Then at $x$ we can compute
\begin{align}
  \Delta u  &= \sum_{i=1}^2 \nb_{e_{i}}\nb_{e_{i}}u \notag\\
                &= \sum_{i=1}^2  \nb_{e_i} \langle \n, e_i\rangle \notag\\
                &= \sum_{i=1}^2  (\langle \bar{\nb}_{{e}_i}  \n, e_i\rangle +  \langle \n, \bar{\nb}_{{e}_i} e_i\rangle)\notag\\
                &= \sum_{i=1}^2 \langle \tanh(u)(e_i - \langle \n, e_i\rangle \n), e_i \rangle + \sum_{i=1}^2 \langle \n, \bar{\nb}_{e_i} e_i \rangle\notag\\
                &=  2\tanh(u) - \tanh(u)(1-\Theta^2) - H\Theta \notag\\
                &= \tanh(u)(1+\Theta^2) - H\Theta,
\end{align}
where we have used Lemma {\ref{generalid}}.
\ep
\br\label{r}
Combining with \eqref{eq-u0}, we have the {\ee} for the hyperbolic distance function $u$ of $S(t)$ along the {\mcf}:
\be\label{eq-r}
u_t - \Delta u = -\tanh(u)(1+\Theta^2)\,,
\ene
which yields similar decay of $u$ as in Theorem \ref{bd-u}.
\er
%%%%%%%%%%%%%%%%%%%%%%%%%%%%%%%
\subsection{Evolution equation for the angle function} In this subsection, we take advantage of the presence of a 
special vector field $V = \cosh(r)\n$ (see \eqref{killing}) to derive the {\ee} for the angle function $\Theta(\cdot,t)$. Recall that, 
on the {\es} $S(t)$, it is given by $\Theta(\cdot,t) = \langle \n,\vnu\rangle (\cdot, t)$. We find the following function more convenient to work with in our hyperbolic setting:
\be
\eta(\cdot,t) = \cosh(u)\Theta(\cdot,t) = \langle V,\vnu\rangle.
\ene
\bl\label{lemma-eta}
The function $\eta(\cdot,t)$ on the {\es} $S(t)$ satisfies the following equation:
\be\label{eq-eta0}
 \Delta \eta = \sinh(u) H -|A|^{2}\eta + \langle V, \nb H\rangle.
\ene
Here $\Delta$ is the Laplace operator on $S(t)$ {\wrt} the induced metric.
\el
\bp
For any point $p \in S(t)$, we choose $\{e_1,e_2\}$ (with $e_3 = \vnu$) to be a local normal frame at $p$.
Then at $p$ we have
     \begin{align*}
     \Delta \eta &= \sum_{i=1}^2 \nb_{e_{i}}\nb_{e_{i}}\langle V,\vnu\rangle\notag\\
                 &= \sum_{i=1}^2 \langle\bnb_{e_{i}}\bnb_{e_{i}}V,\vnu\rangle+
                 2\langle\bnb_{e_{i}}V,\bnb_{e_{i}}\vnu\rangle+
                 \langle V,\bnb_{e_{i}}\bnb_{e_{i}}\vnu\rangle\notag\\
            & = \sum_{i=1}^2 \langle\bnb_{e_{i}}(\sinh(u)e_{i}),\vnu\rangle + 2 \sum_{i,k=1}^2 \sinh(u)\langle e_{i}, a_{ik}e_{k}\rangle + \sum_{i,k=1}^2 \langle V, \bnb_{e_{i}}(a_{ik}e_{k})\rangle  \notag\\
            & = -\sinh(u)H + 2\sinh(u)H + \sum_{i,k=1}^2 a_{ik}\langle V, \bnb_{e_{i}}e_{k}\rangle+a_{ik, i}\langle V,e_{k}\rangle  \notag\\
            & = \sinh(u)H -  \sum_{i,k=1}^2 a_{ik}^2\eta + a_{ii, k}\langle V,e_{k}\rangle\\
            &=  \sinh(u) H -|A|^{2}\eta + \langle V, \nb H\rangle\,,
     \end{align*}
where in the second to the last equality we have used the Codazzi equation and the fact that {\Fm} $M^3$ is of constant
curvature $-1$.

\ep

We next compute $\Delta\Theta$ over the surface $S(t)$.
\bl\label{DTheta} With the above notations, we have:
\begin{align}\label{Theta1}
  \Delta \Theta &= \langle\nabla H, \n\rangle -|A|^{2}\Theta+\tanh(u)(1+\Theta^{2})H
  -\F{\Theta(1-\Theta^{2})}{\cosh^{2}(u)} \notag \\
       &\ \ \ \ -2\tanh(u)\langle\nabla \Theta, \n \rangle -2\tanh^{2}(u)\Theta.
\end{align}
 \el
 \bp
A direct calculation yields
\begin{align}
\Delta \eta &= \Delta (\cosh(u)\Theta) \notag\\
&= \cosh(u)\Delta\Theta + 2\sinh(u)\langle\nb\Theta, \n\rangle +\Theta(\sinh(u)\Delta u+\cosh(u)|\nb u|^2).
\end{align}
Isolating $\Delta\Theta$, we have
  \begin{align*}
  \Delta \Theta &= \F{\Delta\eta}{\cosh(u)}-2\tanh(u)\langle\nb\Theta,\n\rangle-\Theta \tanh(u)\Delta u -\Theta(1-\Theta^{2})\\
  &=\tanh(u)H-|A|^{2}\Theta +\langle \nb H, \n\rangle -2\tanh(u)\langle\nb\Theta, \n\rangle\\
  &\ \ \ \ -\tanh^{2}(u)\Theta(1+\Theta^{2})+\tanh(u)H\Theta^{2}-\Theta(1-\Theta^{2}).
  \end{align*}
Here we have used the fact that $|\nb u|^{2}=1-\Theta^{2}$. Now standard hyperbolic trigonometric identities give
 \beq
 -\Theta(1-\Theta^{2})-\tanh^{2}(u)\Theta(1+\Theta^{2})=-\F{\Theta(1-\Theta^{2})}{\cosh^{2}(u)}-2\tanh^{2}(u)\Theta.
 \eeq
 This completes the proof.
 \ep
 
We can now derive the {\ee} for $\Theta$ along the flow.
\bt\label{ee-theta} With the above notations, we have:
 \begin{align}\label{evo-Theta}
\ppl{\Theta}{t}  - \Delta\Theta &=|A|^{2}\Theta -2\tanh(u) H +2\tanh(u)\langle \nb\Theta, \n\rangle \notag\\
&\ \ \ \ +\F{\Theta(1-\Theta^{2})}{\cosh^{2}(u)}+2\tanh^{2}(u)\Theta.
\end{align}
\et
 \bp
Recall that for the {\mcf} we have $\ppl{}{t}\,\vnu=\nabla H$, and geometrically, one can view $\frac{\p}{\p t}$ as the 
spatial covariant derivative $-H\vnu$ here. Therefore
\beq
   \ppl{\Theta}{t} = \ppl{}{t}\langle\n,\vnu\rangle =\left \langle \ppl{}{t}\vnu, \n \right\rangle + \langle\vnu, \bnb_{-H\vnu}\n\rangle  = \langle\nb H, \n\rangle -H \langle\vnu, \bnb_{\vnu}\n\rangle.
\eeq
Using \eqref{id} we have
\beq
\bnb_{\vnu}\n = \tanh(u)(\vnu - \Theta\n),
\eeq
and thus
\be\label{theta-t}
\ppl{\Theta}{t} = \langle\nb H, \n\rangle -H\tanh(u)(1-\Theta^2).
\ene
Now the conclusion follows from Lemma ~\ref{DTheta}.
\ep

%%%%%%%%%%%%%%%%%%%%%%%%%%%%%%%%%%%%
\subsection{Preserving graphical property} 

In this subsection, we discuss the preservation of the graphical property. We formulate the a priori estimate on the angle $\Theta$ as follows.

\bt\label{main3}
Let $M^3$ be a Fuchsian manifold and $S_0$ be a smooth closed surface which is a geodesic graph over the unique {\tg} surface $\Sigma$ in $M^3$, and suppose there is a positive constant $a_0$ such that $S_0$ lies entirely between $\Sigma(\pm a_0)$. Then whenever the {\ins} $S_0$ satisfies $\Theta_0 \ge \tanh(a_0)$, the {\mcf} with initial surface $S_0$ remains as geodesic graph over $\Sigma$, namely $\Theta(\cdot, t) > 0$ as long as the flow exists.
\et
\bp
We have $-a_0\le u(x, 0) \le a_0$ for any $x\in S_0$. By Theorem \ref{bd-u}, we have:
\be
|u(x, t)| \le \sinh^{-1}(e^{-2t}\sinh(a_0)).
\ene 

We want to find a positive lower bound for $\Theta$ at the initial time to guarantee a positive lower bound for all time and hence the convergence. 

It is equivalent to work with the function $\alpha = \Theta^2$ whose {\ee} is slightly more convenient to work with. With the evolution equation \eqref{evo-Theta} for $\Theta$, we easily deduce the evolution equation for $\alpha$: 
\begin{align}\label{evo-alpha}
\ppl{\alpha}{t}  - \Delta\alpha &=2|A|^{2}\alpha -4\tanh(u)H\Theta +4\tanh(u)\Theta\langle \nb\Theta, \n\rangle \notag\\
&\ \ \ \ +\F{2\alpha(1-\alpha)}{\cosh^{2}(u)}+4\tanh^{2}(u)\alpha - 2|\nb\Theta|^2.
\end{align}

Let $$\phi(t)=\min_{S(t)}\alpha $$ and we only need to consider the case of $\phi\in (0, 1)$ in search of a priori estimate. At the spatial minimum point of $\alpha$, we have $\nabla \Theta=0$ and $\Delta\alpha\ge 0$, and so for $t>0$ (using Hamilton's trick):
\begin{align}
 \F{d\phi}{dt} & \ge \ppl{\alpha}{t}  - \Delta\alpha \notag\\
   &=2|A|^{2}\alpha -4\tanh(u) H \Theta  +\F{2\alpha(1-\alpha)}{\cosh^{2}(u)}+4\tanh^{2}(u)\alpha  \notag\\
     &\ge 2|A|^{2}\Theta^2 -4\sqrt{2}\tanh(|u|)|A|\Theta+\F{2\alpha(1-\alpha)}{\cosh^{2}(u)}+4\tanh^{2}(u)\alpha \notag\\
     &= 2(|A|\Theta-\sqrt{2}\tanh(|u|))^2 -4\tanh^2(u) +\F{2\alpha(1-\alpha)}{\cosh^{2}(u)}+4\tanh^{2}(u)\alpha\notag\\
     &\ge -4(1-\alpha)\tanh^2(u)=-4(1-\phi)\tanh^2(u). \notag
\end{align}

Combining with $|u(x, t)| \le \sinh^{-1}(e^{-2t}\sinh(a_0))$, we end up with the following ordinary differential inequality
\be\label{ivp}
   \left\{
   \begin{aligned}
     \left(\frac{1}{1-\phi}\right)\ddl{\phi}{t}&\ge \frac{-4e^{-4t}\sinh^2(a_0)}{1+e^{-4t}\sinh^2(a_0)}\ \\
     \phi(0)&=\phi_0 \,,
   \end{aligned}
   \right.
\ene
where by assumption we have 
\be\label{initial}
1>\phi_0=\min_{S_0}\alpha\ge \tanh^2(a_0)=\frac{\sinh^2(a_0)}{1+\sinh^2(a_0)}.
\ene
With a proper choice of $\epsilon\in [0, 1)$, we have 
\be\label{phi_0}
\phi_0 =\frac{\epsilon+\sinh^2(a_0)}{1+\sinh^2(a_0)},
\ene 
where $\epsilon>0$ if (\ref{initial}) is a strict inequality.

Straightforward calculations yield:
\be\label{phi-explicit}
\phi(t) \ge \frac{\epsilon+\sinh^2(a_0)e^{-4t}}{1+\sinh^2(a_0)e^{-4t}}> \epsilon
\ene
for all $t\geq 0$ and so we have the following a priori estimate for $t\in [0, \infty):$
\be
\Theta^2(\cdot,t) > \epsilon\ge 0,
\ene
which provides a positive lower bound for the angle $\Theta$ in any finite time interval, and the evolving surface remains as geodesic graph over $\Sigma$ as long as it exists. This completes our proof.
\ep

%%%%%%%%%%%%%%%%%%%%%%
\subsection{Long-time existence and convergence}
Now we can re-assemble the ingredients and complete the proof of Theorem ~\ref{main}.
\bp (of Theorem~\ref{main})
We have shown in Theorem ~\ref{main3} that the {\mcf} \eqref{mcf} stays graphical as long as it exists. This provides the gradient estimate for the {\mcf} for any finite time interval. By the classical theory of parabolic equations in divergent form (for instance \cite{LSU68}), the higher regularity and a priori estimates of the solution follow immediately. This yields the long time existence of the flow by Huisken (\cite{Hui86}). Then by Theorem \ref{bd-u} and the avoidance principle, the continuous convergence of the flow also follows. When the inequality in the assumption is strict, the proof of Theorem ~\ref{main3} gives the uniform estimate of the angle for all time and so the higher order estimates are also uniform for all time, which provides the smooth convergence. The proof of Theorem ~\ref{main} is completed.
\ep

\section{Remarks on the General Case} \label{sec4}

In this section, we discuss the general situation after removing the assumption on the angle function $\Theta_0$ of the initial graph in Theorem \ref{main} and illustrate the relation with possible formation of singularities. In light of the evolution equation (\ref{evo-alpha}), in order to rule out singularities, it will be enough to get a proper bound of $H$ at the point where $\Theta$ takes the spatial minimum in $(0, 1)$, which motivates the following consideration. For the angle function $\Theta$ of any surface in $M^3$, we have the following calculation of its gradient over the surface in general: 
\begin{align}\label{Theta-gradient}
\nabla\Theta
&= \nabla\langle{\n},{\vnu}\rangle \notag\\
&= \langle\bar\nabla_{e_i}\vnu, \n\rangle e_i+\langle\bar\nabla_{e_i} \n, \vnu\rangle e_i \notag\\
&= \langle a_{ij}e_j, \n\rangle e_i+\langle\tanh(u)(e_i-\langle e_i, \n\rangle \n), \vnu\rangle e_i \notag\\
&= a_{ij}\langle e_j, \n\rangle e_i-\tanh(u)\langle e_i, \n\rangle\Theta e_i
\end{align}
where $\{e_1, e_2\}$ is any orthonormal frame at the point of interest on the surface, and we have used Lemma \ref{generalid} in the second to the last step. By choosing $e_1$ and $e_2$ to be the two principal directions so that the second fundamental form is 
$$(a_{ij})={\rm diag}\{a, b\},$$ 
we have 
\begin{align}\label{gradient-Theta}
\nabla\Theta=(a-\tanh(u)\Theta)\langle e_1, \n\rangle e_1+(b-\tanh(u)\Theta)\langle e_2, \n\rangle e_2.
\end{align} 

Notice that the spatial maximum of $\Theta$ is clearly $1$ and obtained at the points with extremal height, where both $e_i$'s are perpendicular to $\n$. Also recall $|\nabla u|^2=\langle e_1, \n\rangle^2+\langle e_2, \n\rangle^2=1-\Theta^2$. Meanwhile, at the spatial minimal point of $\Theta$, denoted by $\theta$ and assumed to be in $(0, 1)$, we see below that exactly one of $e_1$ and $e_2$ is perpendicular to $\n$. Pick a tangent vector of the surface $S$, $\hat e_1$ perpendicular to $\n$, which is unique up to sign as it is also perpendicular to $\vnu$. Then take $\hat e_2$ accordingly so that they form an orthonormal frame of the tangent space of $S$. Use $\hat a_{ij}$ to denote the coefficients for second fundamental form with respect to this basis. Applying (\ref{Theta-gradient}), we have
\beq
\hat e_1(\Theta)=a_{12}\langle \hat e_2, \n\rangle,
\eeq
which implies $a_{12}=0$ since $\langle \hat e_2, \n\rangle^2=1-\theta^2>0$. So such chosen $\hat e_1$ and $\hat e_2$ are also principal directions and can be taken as $e_1$ and $e_2$ as above. So we can have $\langle e_1, \n\rangle=0$. By (\ref{gradient-Theta}), we have $b=\tanh(L)\theta$ where $L$ is the height of the spatial minimal point under consideration.  

Since the flow starts with a graph, singularities can occur only when there is no longer positive lower bound for $\Theta$ by the discussion at the end of Section \ref{sec3}. Thus by defining 
$$T=\sup\left\{t\in (0, \infty)\,|\,\theta\geq C ~\text{in}~[0, t] ~\text{for some}~ C>0\right\}\in (0, \infty],$$  
we know the flow exists as graph exactly in $[0, T)$. In the following, we focus on the case of $T<\infty$, i.e. the flow fails to be graphical in finite time.

In this case, we can choose a time sequence $\{t_i\}_{i=1}^\infty$ approaching $T$ as $i\to\infty$, such that $\theta(t_i)\to 0$ as $i\to\infty$, i.e. $\langle \n, \vnu\rangle\to 0$ at the spatial minimal point $p_i$ for $\Theta$. Now we analyze $p_i\in S(t_i)$ at time $t_i$ more carefully. For simplicity of notations, we frequently omit the index $i$ below, and the limit is always taken as $i\to\infty$. We already know that there is one principal direction $e_1$ such that $\langle e_1, \n\rangle=0$ at the spatial minimal point. Geometrically, $e_1$ is the direction of the curve as the intersection of the evolving graph $S(t)$ and the equidistance graph $\Sigma(L)$ where $L$ is the height of the spatial minimal point. As $\langle e_1, \n\rangle^2+\langle e_2, \n\rangle^2=1-\theta^2$, we have $\langle e_2, \n\rangle^2\to 1$, i.e. $e_2\to \n$ by reversing $e_2$ if necessary. Consider the geodesics on $S(t_i)$ starting at $p_i$ in the direction of $e_2$. Since at $p_i$ we have $\langle \bar\nabla_{e_2}e_2, \vnu\rangle=-\langle \bar\nabla_{e_2}\vnu, e_2\rangle=-b=-\tanh(L)\theta$, which approaches $0$ by the decay of height and $\theta(t_i)\to 0$, we have $\bar\nabla_{e_2}e_2\to 0$ where by abuse of notation, $e_2$ also stands for the unit tangent vector field along the geodesic. Together with $e_2\to \n$, we know that in the infinitesimal way at $p_i$, this geodesic on $S(t_i)$ approaches the $r$-curve which is geodesic of $M^3$ in the direction $\n$. Intuitively, this is the consequence of the loss of graphical property.      

Meanwhile, the ``reason" for the loss of graphical property should be the ``relative" blow-up of the principle curvature in the $e_1$ direction, namely the quantity $a$ in (\ref{gradient-Theta}). We hope to illustrate this point in the following. Using (\ref{evo-Theta}) and $b=\tanh(L)\theta$, we have
\be\label{general-evo-theta}
\frac{d\theta}{d t}\geq \bigl(a^2+\tanh^2(L)\theta^2\bigr)\theta-2a\tanh(L)+\frac{\theta(1-\theta^2)}{\cosh^2(L)}
\ene
If $|a|\leq-\theta\log\theta+C\theta$ for some $C>0$ at any spatial minimal point of $\Theta$ as long as the flow exists, then we have  
\begin{align}\label{example-evo-theta}
\frac{d\theta}{d t}\geq C\theta\log\theta-C\theta.
\end{align}
Direct calculation yields that $\theta\geq Ce^{-Ce^{Ct}}>0$, which rules out the loss of graphical property in any finite time and we have the long time existence together with the continuous convergence. Motivated by this, for the case of $T<\infty$ and some $C>0$, we set 
\be
\I=\{t\in [0, T)\,|\,|a|>-\theta\log\theta+C\theta~\text{at all spatial minimal points of}~\Theta(\cdot, t)\}.
\ene
Note that $\I$ has the closure in $\R$ containing $T$, since otherwise we can derive a positive lower bound for $\theta$ for any finite time interval as above, contradicting $T<\infty$. So we can choose the time sequence $\{t_i\}_{i=1}^\infty\subset \I$ approaching $T$ as $i\to \infty$ and consider at the point $q_i$ where $\Theta$ achieves the spatial minimum on $S(t_i)$ and $|a|>-\theta\log\theta+C\theta$. We can further make sure that $\theta\to 0$ for this time sequence, since otherwise $\theta$ will have a uniform positive lower bound for $t\in \I$ close to $T$, and $\theta$ will then have a uniform positive lower bound  for $t\in [0, T)\setminus \I$  by applying the above argument for the interior of $[0, T)\setminus \I$, contradicting the choice of $T$. Here we make use of the continuity of $\theta$ with respect to time.

Since $\theta\to 0$, the scale of $b=\tanh(L)\theta$ is small and way smaller than that of $a$, i.e. ``relative" blow-up. In other words, after a proper ``blowing-down" (by the scale of $\theta(-\log\theta)^{1/2}$, for example), we have the violation of graphical property modelled as a cylinder with the circle on the equidistant surface $\Sigma(L)$ and pointing in the $\n$ direction in the infinitesimal way. 

By the discussion at the end of Section \ref{sec3}, the blow-up of $|A|$, i.e. formation of flow singularities as surface, can't occur before the degeneration of $\Theta$ and certainly might not happen at the same place. This adds to the intriguing features about the singularities. In future works, we hope to provide more precise local and global understanding for such singularities, aiming at either ruling them out or obtaining interesting examples.  
%%%%%%%%%%%%%%%%%%%%

\section{Appendix}

In this appendix, we give a detailed proof for the {\ee} for the angle function in Theorem ~\ref{old-eq}, i.e. the equations \eqref{ev-theta} and \eqref{ev-theta-1} for our Riemannian setting. 
The calculation is carried out for the {\mcf} of graphical hypersurfaces of general dimension $n$, in the ambient manifold $M^{n+1}$ with a general warped product metric.

We still use $\n$ and $\vnu$ to denote the unit normal vectors for the warped product foliation and the evolving hypersurface.  Also we sum over all repeated indices in this section.

The following computation is done for $F(p, t)$ for time $t$. We choose the normal frame $\{e_i\}_{i=1}^n$ for the evolving hypersurface. Then we require 
$L_\n e_i=[\n, e_i]=0$ to extend the frame to a neighborhood of $F(p, t)$ in $M^{n+1}$, i.e. using $\n$ to generate a family of hypersurfaces with the initial one being the 
evolving hypersurface at time $t$. Then the vector field $\vnu$ below means the normal vector field for this family of hypersurfaces. This won't affect the result for $\Delta\Theta$ at 
$F(p, t)$ for time $t$. 
\begin{align*}
\Delta\Theta &= e_ie_i\<\vnu, \n\> \\
&= e_i\(\<\vnu, \bar\nabla_{e_i}\n\>+\<\bar\nabla_{e_i}\vnu, \n\>\) \\
&= e_i\(\<\vnu, \bar\nabla_\n e_i\>+\<\bar\nabla_{e_i}\vnu, e_j\>\cdot\<\n, e_j\>\) \\
&= e_i\(\<\vnu, \bar\nabla_\n e_i\>+\<\bar\nabla_{e_j}\vnu, e_i\>\cdot\<\n, e_j\>\) \\
&= \<\bar\nabla_{e_i}\vnu, \bar\nabla_\n e_i\>+\<\vnu, \bar\nabla_{e_i}\bar\nabla_\n e_i\>+\<\bar\nabla_{e_i}\bar\nabla_{e_j}\vnu, e_i\>\cdot\<\n, e_j\> \\
&+\<\bar\nabla_{e_j}\vnu, \bar\nabla_{e_i}e_i\>\cdot\<\n, e_j\> +\<\bar\nabla_{e_j}\vnu, e_i\>\cdot\<\bar\nabla_{e_i}\n, e_j\>
+\<\bar\nabla_{e_j}\vnu, e_i\>\cdot\<\n, \bar\nabla_{e_i}e_j\> \\
&= \<\bar\nabla_{e_i}\vnu, \bar\nabla_\n e_i\>+\<\vnu, \bar\nabla_{e_i}\bar\nabla_\n e_i\>+\<\bar\nabla_{e_i}\bar
\nabla_{e_j}\vnu, e_i\>\cdot\<\n, e_j\> \\
&\ \ +\<\bar\nabla_{e_j}\vnu, e_i\>\cdot\<\bar\nabla_{e_i}\n, e_j\>+\<\bar\nabla_{e_j}\vnu, e_i\>\cdot\<\n, \bar
\nabla_{e_i}e_j\>, 
\end{align*}
where we have used $\bar\nabla_{e_i} \n=\bar\nabla_\n e_i$ for the third equality, $a_{ij}=\<\bar\nabla_{e_i}\vnu, 
e_j\>=\<\bar\nabla_{e_j}\vnu, e_i\>$ for the fourth equality and $\bar\nabla_{e_i}e_i=-a_{ii}\vnu=-H\vnu$ for the last 
equality. For these terms, we have
\beq
\<\vnu, \bar\nabla_{e_i}\bar\nabla_\n e_i\>=\<R(e_i, \n)e_i, \vnu\>+\<\vnu, \bar\nabla_\n\bar\nabla_{e_i} e_i\>,
\eeq
where $R$ is the Riemannian curvature tensor, and we also have 
\begin{equation}
\begin{split}
e_j(H)
&=e_j\<\bar\nabla_{e_i}\vnu, e_i\> \\
&= \<\bar\nabla_{e_j}\bar\nabla_{e_i}\vnu, e_i\>+\<\bar\nabla_{e_i}\vnu, \bar\nabla_{e_j}e_i\> \\
&= -\<R(e_i, e_j)\vnu, e_i\>+\<\bar\nabla_{e_i}\bar\nabla_{e_j}\vnu, e_i\>, \nonumber
\end{split}
\end{equation}
using $[\n, e_i]=0$, $\bar\nabla_{e_i}e_j=\bar\nabla_{e_j}e_i=-a_{ij}\vnu$ and $[e_i, e_j]=0$. We also find:
\begin{equation} 
\begin{split}
\<\bar\nabla_{e_i}\vnu, \bar\nabla_\n e_i\>
&= \<\bar\nabla_{e_i}\vnu, e_j\>\cdot\<\bar\nabla_\n e_i, e_j\> \\
&= \<\bar\nabla_{e_j}\vnu, e_i\>\cdot\<\bar\nabla_{e_i}\n, e_j\> \\
&= \<\bar\nabla_\n e_i, e_j\>a_{ij} \\
&= \frac{a_{ij}}{2}\n\(\<e_i, e_j\>\). \nonumber
\end{split}
\end{equation}
So the previous computation for $\Delta\Theta$ can be continued as follows:
\begin{equation}
\begin{split}
\Delta\Theta
&= \<\bar\nabla_{e_i}\vnu, \bar\nabla_\n e_i\>+\<\vnu, \bar\nabla_{e_i}\bar\nabla_\n e_i\>+\<\bar\nabla_{e_i}\bar
\nabla_{e_j}\vnu, e_i\>\cdot\<\n, e_j\> \\
&~~~~~~~~~ +\<\bar\nabla_{e_j}\vnu, e_i\>\cdot\<\bar\nabla_{e_i}\n, e_j\>+\<\bar\nabla_{e_j}\vnu, e_i\>\cdot\<\n, \bar
\nabla_{e_i}e_j\> \\
&= a_{ij}\n\(\<e_i, e_j\>\)+\<R(e_i, \n)e_i, \vnu\>+\<\vnu, \bar\nabla_\n\bar\nabla_{e_i} e_i\> \\
&~~~~~~ +\<\n, e_j\>\cdot\(e_j(H)+\<R(e_i, e_j)\vnu, e_i\> \)+\<\bar\nabla_{e_j}\vnu, e_i\>\cdot\<\n, \bar\nabla_{e_i}e_j\>. \nonumber
\end{split}
\end{equation}
We consider each term separately below:
$$a_{ij}\n\(\<e_i, e_j\>\)=-a_{ij}\n(g^{ij}),$$
$$\<R(e_i, \n)e_i, \vnu\>=-\Ric(\n, \vnu),$$
$$\<\vnu, \bar\nabla_\n\bar\nabla_{e_i} e_i\>=\n\(\<\vnu, \bar\nabla_{e_i} e_i\>\),$$
$$\<\n, e_j\>\cdot e_j(H)=\<\n, \nabla H\>,$$
$$\<\n, e_j\>\cdot\<R(e_i, e_j)\vnu, e_i\>=\Ric(\n^{\l}, \vnu),$$
$$\<\bar\nabla_{e_j}\vnu, e_i\>\cdot\<\n, \bar\nabla_{e_i}e_j\>=a_{ij}\cdot\<\n, -a_{ij}\vnu\>=-\sum_{i, j}|a_{ij}|^2\cdot\<\n, \vnu\>,$$
where $\bar\nabla_{e_i}e_i=-a_{ii}\vnu=-H\vnu$ is used for the third one, $(g^{ij})$ is the inverse matrix of $(g_{ij}=\<e_i, e_j\>)$ and 
$\n^{\l}$ is the projection of $\n$ in the direction of the evolving hypersurface.

\br The equality $\<\vnu, \bar\nabla_{e_i} e_i\>=-H$ holds only at $F(p, t)$ for time $t$, and so $\n\(\<\vnu, \bar\nabla_{e_i} e_i\>\)$ is 
NOT equal to $-\n(H)$. Nevertheless, we still have $\<\vnu, \bar\nabla_{e_i} e_i\>=-\<\bar\nabla_{e_i}\vnu, e_i\>$ by the construction 
of $\vnu$ at the beginning of this appendix.
\er

Now we can finish the computation for $\Delta\Theta$:
\begin{equation}
\begin{split}
\Delta\Theta
&= -a_{ij}\n(g^{ij})+\n\(\<\vnu, \bar\nabla_{e_i}e_i\>\)-\<\n, \vnu\>\cdot\Ric(\vnu, \vnu)\\
& \ \ \ \ \ \ \ +\<\n, \nabla H\>-\sum_{i, j}|a_{ij}|^2\cdot\<\n, \vnu\> \\
&= -\n(H_\n)+(\Ric(\vnu, \vnu)-\sum_{i, j}|a_{ij}|^2)\Theta+\<\n, \nabla H\>, 
\nonumber
\end{split}
\end{equation}
where $\n(H_\n)=a_{ij}\n(g^{ij})-\n(\<\vnu, \bar\nabla_{e_i} e_i\>)$.

Next we further clarify the term $\n(H_\n)$. It stands for the variation of mean curvature for the family of hypersurfaces starting with the evolving 
hypersurface under consideration at time $t$ and flowing out by the vector field $\n$. This is the same family of hypersurfaces considered in 
the previous calculation, and we are only interested in the initial hypersurface which is the hypersurface evolving along the {\mcf} at time $t$. 
Clearly, we have 
\beq
\n(H_\n)=\n(g^{ij}a_{ij})=g^{ij}\n(a_{ij})+a_{ij}\n(g^{ij})=\n(a_{ii})+a_{ij}\n(g^{ij})
\eeq
at the point since $(g^{ij})$ is the identity matrix at the point under consideration, and $a_{ii}$ is not equal to $H$ nearby. 
The computation for $\n(H_\n)$ is as follows, still for just that point.
\begin{equation}
\begin{split}
\n(H_\n)
&= \n(a_{ii})+a_{ij}\n(g^{ij}) \\
&= \n\(\<\bar\nabla_{e_i}\vnu, e_i\>\)-a_{ij}\n\(\<e_i, e_j\>\) \\
&= \<\bar\nabla_\n\bar\nabla_{e_i}\vnu, e_i\>+\<\bar\nabla_{e_i}\vnu, \bar\nabla_\n e_i\>-a_{ij}L_\n g(e_i, e_j) \\
&= -a_{ij}(L_\n g)(e_i, e_j)+e_i\(\<\vnu, \bar\nabla_\n e_i\>+\<e_i, \bar\nabla_\n \vnu\>\)-\<\bar\nabla_{e_i}e_i, 
\bar\nabla_\n\vnu\> \\
&~~~~ -\<\vnu, \bar\nabla_{e_i}\bar\nabla_\n e_i\>-\<e_i, \bar\nabla_{e_i}\bar\nabla_\n\vnu\>+\<e_i,\bar\nabla_\n\bar
\nabla_{e_i}\vnu\> \\
&=  -a_{ij}(L_\n g)(e_i, e_j)-\<\vnu, \bar\nabla_{e_i}\bar\nabla_\n e_i\>+\<R(\n, e_i)\vnu, e_i\> \\
&= -a_{ij}(L_\n g)(e_i, e_j)-e_i\(\<\vnu, \bar\nabla_\n e_i\>\)+\<\bar\nabla_{e_i}\vnu, \bar\nabla_\n e_i\>+\<R(\vnu, e_i)\n, e_i\> \\
&= -a_{ij}(L_\n g)(e_i, e_j)-e_i\(\<\vnu, \bar\nabla_\n e_i\>\)+\<\bar\nabla_{e_i}\vnu, \bar\nabla_\n e_i\> \\
&~~~~ +\<\bar\nabla_{\vnu}\bar\nabla_{e_i}\n, e_i\>-\<\bar\nabla_{e_i}\bar\nabla_{\vnu}\n, e_i\>-\<\bar\nabla_{\bar\nabla_{\vnu}e_i}\n, e_i\>
+\<\bar\nabla_{\bar\nabla_{e_i}\vnu}\n, e_i\>, \nonumber
\end{split}
\end{equation}
where $\n(g^{ij})=-g^{ik}\n(g_{k\ell})g^{\ell j}$ and $(g^{ij})$ being identity at the point are used for the first equality; $[\n, e_i]=0$ is used for the 
third equality; $\<e_i, \vnu\>=0$, $|\vnu|=1$ and at $F(p, t)$, $\bar\nabla_{e_i}e_i=-a_{ii}\vnu=-H\vnu$ are used for the fifth equality. 

We have a few more terms to sort out.
\begin{equation}
\begin{split}
-a_{ij}L_\n g(e_i, e_j)
&= -\<e_i, \bar\nabla_{e_j}\vnu\>\(\<\bar\nabla_{e_i}\n, e_j\>+\<\bar\nabla_{e_j}\n, e_i\>\) \\
&= -\<\bar\nabla_{\bar\nabla_{e_j}\vnu}\n, e_j\>-\<\bar\nabla_{e_j}\n, \bar\nabla_{e_j}\vnu\>, \nonumber 
\end{split}
\end{equation}
$$\frac{1}{2}\vnu\(L_\n g(e_i, e_i)\)= \vnu\(\<\bar\nabla_{e_i}\n, e_i\>\)= \<\bar\nabla_{\vnu}\bar\nabla_{e_i}\n, e_i\>+\<\bar\nabla_{e_i}\n, \bar\nabla_{\vnu}e_i\>,$$
$$-L_\n g(\bar\nabla_{\vnu}e_i, e_i)=-\<\bar\nabla_{\bar\nabla_{\vnu}e_i}\n, e_i\>-\<\bar\nabla_{e_i}\n, \bar\nabla_
{\vnu}e_i\>,$$
\begin{equation}
\begin{split}
-e_i\(L_\n g(\vnu, e_i)\) 
&= -e_i\(\<\bar\nabla_{\vnu}\n, e_i\>+\<\bar\nabla_{e_i}\n, \vnu\>\) \\
&= -\<\bar\nabla_{e_i}\bar\nabla_{\vnu}\n, e_i\>-\<\bar\nabla_{\vnu}\n, \bar\nabla_{e_i}e_i\>-e_i\(\<\bar\nabla_\n e_i, \vnu\>\). \nonumber 
\end{split}
\end{equation}
Now it is easy to calculate: 
\begin{equation}
\begin{split}
\n(H_\n)
&= \frac{1}{2}\vnu\(L_\n g(e_i, e_i)\)-L_\n g(\bar\nabla_{\vnu}e_i, e_i)-e_i\(L_\n g(\vnu, e_i)\)+\<\bar\nabla_{\vnu}\n, \bar\nabla_{e_i}e_i\> \\
&= \frac{1}{2}(\bar\nabla_{\vnu}L_\n g)(e_i, e_i)-e_i\(L_\n g(\vnu, e_i)\)+\<\bar\nabla_{\vnu}\n, \bar\nabla_{e_i}e_i\> \\
&= \frac{1}{2}(\bar\nabla_{\vnu}L_\n g)(e_i, e_i)+\<\bar\nabla_{\vnu}\n, \bar\nabla_{e_i}e_i\> \\
&\ \ \ \ \ \   -(\bar\nabla_{e_i}L_\n g)(\vnu, e_i)-L_\n g(\bar\nabla_{e_i}\vnu, e_i)-L_\n g(\vnu, \bar\nabla_{e_i}e_i) \\
&= \frac{1}{2}(\bar\nabla_{\vnu}L_\n g)(e_i, e_i)+\<\bar\nabla_{\vnu}\n, \bar\nabla_{e_i}e_i\> \\
&\ \ \ \ \ \  -(\bar\nabla_{e_i}L_\n g)(\vnu, e_i)-L_\n g(a_{ij}e_j, e_i)-L_\n g(\vnu, -H\vnu). \nonumber
\end{split}
\end{equation}
In light of $\<\bar\nabla_{\vnu}\n, \bar\nabla_{e_i}e_i\>=-H\<\bar\nabla_{\vnu}\n, \vnu\>=-\frac{1}{2}HL_\n g(\vnu, \vnu)$, 
we conclude
$$\n(H_\n)=\frac{1}{2}(\bar\nabla_{\vnu}L_\n g)(e_i, e_i)-(\bar\nabla_{e_i}L_\n g)(\vnu, e_i)-a_{ij}L_\n g(e_i, e_j)+\frac{1}{2}HL_\n g(\vnu, \vnu).$$
Noticing $\<\vnu, \bar\nabla_{\vnu} \n\>=\frac{1}{2}L_\n g(\vnu, \vnu)$, we have 
\beq
\n(H_\n)=\frac{1}{2}(\bar\nabla_{\vnu}L_\n g)(e_i, e_i)-(\bar\nabla_{e_i}L_\n g)(\vnu, e_i)-a_{ij}L_\n g(e_i, e_j)+H\<\vnu, \bar\nabla_{\vnu} \n\>.
\eeq

We note that the advantage of computing this way is that the terms now depend mostly on the evolving hypersurface. The vector field $\n$ only appears in $L_\n g$. In the following, we compute $\frac{\p \Theta}{\p t}$ in detail. We use the coordinate system used in \cite{Hui86}, i.e. a normal coordinate system 
$\{y_\alpha\}$ for $F(p, t)$ in $M$ with the frame vector for the first coordinate is $-\vnu$ at time $t$.

Let $\vnu=\nu^\alpha\frac{\p}{\p y^\alpha}$ and $\n=n^\alpha\frac{\p}{\p y^\alpha}$. We have $\frac{\p \Theta}{\p t}=\frac{\p (g^{\alpha\beta}\nu^\alpha n^\beta)}{\p t}$. 
There are three terms from Leibniz rule. 
\beq
\frac{\p g_{\alpha\beta}}{\p t}= \frac{\p}{\p t}\(\left\<\frac{\p}{\p y^\alpha}, \frac{\p}{\p y^{\beta}}\right\>\)= -H\vnu\(\left\<\frac{\p}{\p y^\alpha}, \frac{\p}{\p y^{\beta}}\right\>\)= 0,
\eeq
because the Christoffel symbols vanish at the point. Define $\frac{\p \vnu}{\p t}$ to be $\frac{\p \nu^\alpha}{\p t}\frac{\p}{\p y^\alpha}$ and $\frac{\p \n}{\p t}$ to be $\frac{\p n^\alpha}{\p t}\frac{\p}{\p y^\alpha}$, and we see 
\beq
\frac{\p \Theta}{\p t}=\<\frac{\p \vnu}{\p t}, \n\>+\<\vnu, \frac{\p \n}{\p t}\>.
\eeq 
We have $\frac{\p \vnu}{\p t} = \nabla H$, and 

%Since $\frac{\p g_{\alpha\beta}}{\p t}=0$ and $\<\vnu, \vnu\>=1$, one easily sees $\<\frac{\p \vnu}{\p t}, \vnu\>=0$, and so $\frac{\p \vnu}{\p t}$ is a 
%tangent vector for the evolving hypersurface. One can then choose a coordinate system for the reference hypersurface around $p$, $\{x^i\}$. We also have 
%$$\frac{\p g_{\alpha\beta}}{\p x^k}=\frac{\p F}{\p x^k}\(\<\frac{\p}{\p y^\alpha}, \frac{\p}{\p y^{\beta}}\>\)=0$$
%by the choice of $\{y^\alpha\}$. Now one computes $\frac{\p \vnu}{\p t}$ as follows:
%\begin{equation}
%\begin{split}
%\frac{\p \vnu}{\p t}
%&= g^{k\ell}\<\frac{\p \vnu}{\p t}, \frac{\p F}{\p x^k}\>\frac{\p F}{\p x^\ell} \\
%&= -g^{k\ell}\<\vnu, \frac{\p^2 F}{\p t\p x^k}\>\frac{\p F}{\p x^\ell} \\
%&= -g^{k\ell}\<\vnu, \frac{\p (-H\nu^\alpha)}{\p x^k}\frac{\p}{\p y^\alpha}\>\frac{\p F}{\p x^\ell} \\
%&= g^{k\ell}\frac{\p H}{\p x^k}\frac{\p F}{\p x^\ell} \\
%&= \nabla H, \nonumber
%\end{split}
%\end{equation}
%where $\frac{\p F}{\p x^k}=\frac{\p F^\alpha}{\p x^k}\frac{\p}{\p y^\alpha}$, $\frac{\p F}{\p t}=\frac{\p F^\alpha}{\p t}\frac{\p}{\p y^\alpha}=-H\nu^\alpha\frac{\p}{\p y^\alpha}$ and $\frac{\p^2 F}{\p t\p x^k}=\frac{\p^2 F^\alpha}{\p t\p x^k}\frac{\p}{\p y^\alpha}=\frac{\p (-Hv^\alpha)}{\p x^k}\frac{\p}{\p y^\alpha}$. The second equality above is true because $\<\vnu, \frac{\p F}{\p x^k}\>=0$ and $\frac{\p g_{\alpha\beta}}{\p t}=0$. For the fourth equality, one makes use of $\<\vnu, \vnu\>=1$ and $\frac{\p g_{\alpha\beta}}{\p x^k}=0$.

\beq
\frac{\p \n}{\p t}=\frac{\p n^\alpha}{\p t}\frac{\p}{\p y^\alpha}=-H\vnu(n^\alpha)\frac{\p}{\p y^\alpha}=-H\bar\nabla_{\vnu} \n,
\eeq
where the last equality is true again by the choice of $\{y^\alpha\}$.

Finally, we can conclude (\ref{ev-theta}) and (\ref{ev-theta-1}).

%%%%%%%%%%%%%%%%%%%%
\bibliographystyle{amsalpha}
\bibliography{ref-fuchsian}
\end{document}